\documentclass[12pt]{amsart}

\def\1{{\bf 1}}
\def\2{{\bf 2}}

\def\GG{{\mathcal G}}

\title{Two way subtable sum problems and quadratic Gr\"obner bases}
\author{Hidefumi Ohsugi \and Takayuki Hibi}
\date{}
\address{
Department of Mathematics\\
College of Science\\
Rikkyo University\\
Toshima, Tokyo 171-8501, Japan}
\email{ohsugi@rkmath.rikkyo.ac.jp}

\address{
Department of Pure and Applied Mathematics\\
Graduate School of Information Science and Technology\\
Osaka University \\
Toyonaka, Osaka 560-0043, Japan}
\email{hibi@math.sci.osaka-u.ac.jp}

\begin{document}

\begin{abstract}
Hara, Takemura and Yoshida 
\cite{HTY} discussed toric ideals arising from two way
subtable sum problems and showed that 
these toric ideals are generated by quadratic binomials
if and only if the subtables are either diagonal or triangular.
In the present paper, we show that
if the subtables are either diagonal or triangular,
then their toric ideals possess quadratic Gr\"obner bases. 
\end{abstract}

\maketitle

Fix positive integers $m$ and $n$ and 
$T=\{ (i,j) \ | \ 1 \leq i \leq m , 1 \leq j \leq n  \}$.
Let $K$ be a field and $K[\{u_i\}_{1 \leq i \leq m} \cup
\{v_j\}_{1 \leq j \leq n} \cup \{ w , t \} ]$
be the polynomial ring in $m + n +2$ variables over $K$.
Given a subset $S$ of $T$,
let $R_{S}$ denote the semigroup ring generated by 
those monomials $u_i v_j w$ with $(i,j) \in S$ and
those monomials $u_i v_j t$ with $(i,j) \notin S$.
Let $K[X] = K[\{ x_{i,j} \}_{(i,j) \in T}]$
denote the polynomial ring in $m n$ variables over $K$.
Define the surjective map $\pi: K[X] \rightarrow R_{S}$
with $\pi(x_{i,j}) = u_i v_j w$ if $(i,j) \in S$
and $\pi(x_{i,j}) = u_i v_j t$ if $(i,j) \notin S$.
We call the kernel of $\pi$ 
{\it the toric ideal of two way
subtable sum problems associated with $S$}
and 
denote it by $I_{S}$.
We refer the reader to \cite{CLO} and \cite{Stu}
for fundamental facts on Gr\"obner bases
and toric ideals.

A subset $S \subset T$ is called $2 \times 2$ {\it block diagonal}
if there exist integers $r$, $c$ such that 
$$S=
\{ (i,j) \ | \ 1 \leq i \leq r, 1 \leq j \leq c \}
\cup 
\{ (i,j) \ | \ r < i \leq m, c <  j \leq n \}
$$ after an appropriate interchange of rows and columns.
A subset $S \subset T$ is called {\it triangular}
if $S$ satisfies the condition $(*)$
\begin{center}
$(i,j) \in S \Rightarrow (i',j') \in S$
for all $1 \leq i' \leq i$ and $1 \leq j' \leq j$
\end{center}
after an appropriate interchange of rows and columns.
Hara, Takemura and Yoshida \cite{HTY}
showed that $I_S$ is generated by quadratic binomials 
if and only if $S$ is either $2 \times 2$ block
diagonal or triangular.

Our main theorem is as follows:

\bigskip

\noindent
{\bf Theorem.}
{\it
Work with the same notation as above.
Then the following conditions are equivalent:
\begin{itemize}
\item[(i)]
the toric ideal $I_S$ is generated by quadratic binomials;
\item[(ii)]
the toric ideal $I_S$ possesses a squarefree initial ideal;
\item[(iii)]
the toric ideal $I_S$ possesses a quadratic Gr\"obner basis;
\item[(iv)]
the semigroup ring $R_S$ is normal;
\item[(v)]
the semigroup ring $R_S$ is Koszul;
\item[(vi)]
the subset $S$ is either $2 \times 2$ block diagonal or triangular.
\end{itemize}
}

\bigskip

\noindent
{\it Proof.}
First of all, (i) $\Leftrightarrow $ (vi)
is proved in \cite{HTY}.
Moreover, in \cite{HTY}, it is proved that
if $S$ is neither diagonal nor triangular,
then there exists an indispensable binomial $f^{(+)}  - f^{(-)}\in I_S$
such that neither $f^{(+)} $ nor $f^{(-)}$ is squarefree.
Hence, by \cite[Corollary 6.3]{OHconti},
it follows that (iv) $\Rightarrow $ (vi).
On the other hand, in general, it is known that
(ii) $\Rightarrow $ (iv) and that
(iii) $\Rightarrow $ (v) $\Rightarrow $ (i).
Since $R_S$ is generated by monomials of the same degree,
it follows from the proof of \cite[Proposition 1.6]{quadratic} that
(iii) $\Rightarrow $ (ii).
We now prove that (vi) $\Rightarrow $ (iii).
Suppose (vi) holds.

Suppose that $S \subset T$ is triangular and satisfies $(*)$.
Then $I_S$ is generated by 
$$
{\mathcal G}
=
\{
x_{i,\ell}
x_{j,k}
-
x_{i,k}
x_{j,\ell}
\ 
| 
\ 
1\leq i < j \leq m,
1\leq k < \ell \leq n,
\pi(x_{i,\ell}
x_{j,k}) = \pi(x_{i,k}
x_{j,\ell}
)
\}.
$$
(Note that all monomials appearing in $g \in {\mathcal G}$ are squarefree.)
Throughout the proof, we will perform computations using
Buchberger's criterion for ${\mathcal G}$.
These computations, while routine,
are often very long; as such, we will omit the details and
state only the results.
Since $S$ is triangular,
the following holds:
\begin{itemize}
\item[($\sharp$)]
$x_{i,\ell}
x_{j,k}
-
x_{i,k}
x_{j,\ell}$
($
1\leq i < j \leq m$,
$
1\leq k < \ell \leq n
$)
does not belong to ${\mathcal G}$
if and only if
$S \cap \{(i,\ell),(j,k),(i,k),(j,\ell)\}$ is
either
$\{(i,k)\}$
or
$\{(i,k),(i,\ell),(j,k)\}$.
\end{itemize}
Let $<$ be the lexicographic order on $K[X]$
induced by
$x_{m,1} > x_{m,2} > \cdots > x_{m,n}
> x_{m-1,1} > \cdots > x_{m-1,n} > \cdots > x_{1,1} > \cdots > x_{1,n}
.$
Then the initial monomial of 
$x_{i,\ell}
x_{j,k}
-
x_{i,k}
x_{j,\ell}
$
($
1\leq i < j \leq m$,
$
1\leq k < \ell \leq n
$)
is $x_{i,\ell}
x_{j,k}$.

Suppose that the $S$-polynomial $f$ of 
$g_1, g_2 \in {\mathcal G}$ is not reduced to $0$ by ${\mathcal G}$.
By \cite[Ch.2 \S 9 Proposition 4]{CLO}, the initial monomials of $g_1$ and $g_2$
are not relatively prime.
Hence $f$ is a cubic binomial.
If the monomials of $f$ have a common variable, then 
%
$f$ is reduced to $0$ by ${\mathcal G}$.
Let $f = 
x_{i_1 ,\ell_1} x_{i_2 , \ell_2 } x_{i_3 , \ell_3 }- x_{i_1' ,\ell_1'} x_{i_2' , \ell_2' } x_{i_3' , \ell_3' }$
with the initial monomial $x_{i_1 ,\ell_1} x_{i_2 , \ell_2 } x_{i_3 , \ell_3 }$.
Since $\pi(x_{i_1 ,\ell_1} x_{i_2 , \ell_2 } x_{i_3 , \ell_3 })= \pi
(x_{i_1' ,\ell_1'} x_{i_2' , \ell_2' } x_{i_3' , \ell_3' })$, 
we have $\{i_1,i_2,i_3 \} = \{i_1',i_2',i_3' \}$ and $\{\ell_1,\ell_2,\ell_3\} = \{\ell_1',\ell_2',\ell_3'\} $.
Moreover since the monomials of $f$ have no common variable,
$|\{i_1,i_2,i_3 \}  | = |\{\ell_1,\ell_2,\ell_3\}|=3$.
We may assume $i_1 = i_1' < i_2  = i_2' < i_3 = i_3' $.
By the definition of $<$, 
we have $\ell_3 <  \ell_3' \in \{\ell_1,\ell_2\}$.
For example, if $\ell_1 < \ell_3 < \ell_2$, then $\ell_3' = \ell_2$
and hence, by $\ell_1 \neq \ell_1'$, 
we have $(\ell_1',\ell_2' ) = (\ell_3, \ell_1)$.
(This is the case (6) below.)
By such an argument, $f$ is one of the following:
%
\begin{itemize}
\item[(1)]
$x_{i_1 ,j_ 2} x_{i_2 , j_3 } x_{i_3 , j_1 }- x_{i_1 , j_3 } x_{i_2 , j_1 } x_{i_3 , j_2 }$
\item[(2)]
$x_{i_1 ,j_3 } x_{i_2 , j_2 } x_{i_3 , j_1 }- x_{i_1 , j_2 } x_{i_2 , j_1 } x_{i_3 , j_3}$
\item[(3)]
$x_{i_1 ,j_3 } x_{i_2 , j_2 } x_{i_3 , j_1 }- x_{i_1 , j_1 } x_{i_2 , j_3 } x_{i_3 , j_2 }$
\item[(4)]
$x_{i_1 ,j_3 } x_{i_2 , j_1 } x_{i_3 , j_2 }- x_{i_1 , j_1 } x_{i_2 , j_2 } x_{i_3 , j_3 }$
\item[(5)]
$x_{i_1 ,j_ 2} x_{i_2 , j_3 } x_{i_3 , j_1 }- x_{i_1 , j_1 } x_{i_2 , j_2 } x_{i_3 , j_3 }$
\item[(6)]
$x_{i_1 ,j_1 } x_{i_2 , j_3 } x_{i_3 , j_2 }- x_{i_1 , j_2 } x_{i_2 , j_1 } x_{i_3 , j_3 }$
\end{itemize}
where $1 \leq i_1 < i_2 < i_3 \leq m$ and $1 \leq j_1 < j_2 < j_3 \leq n$.

In each case we will show that $f$ is reduced by ${\mathcal G}$ to a cubic
binomial $h$ where the monomials of $h$ have a common variable.
Since $I_S$ is prime, it will follow that the quadratic factor of $h$
belongs to ${\mathcal G}$, and thus that $f$ is reduced to $0$ by ${\mathcal G}$.

(1)
Suppose that none of
$g_1 = x_{i_1 , j_3 } x_{i_3 , j_2 } -  x_{i_1 , j_2 } x_{i_3 , j_3 }$,
$g_2 = x_{i_1 , j_2 } x_{i_3 , j_1 } -  x_{i_1 , j_1 } x_{i_3 , j_2 }$ and
$g_3 = x_{i_2 , j_3 } x_{i_3 , j_1 } -  x_{i_2 , j_1 } x_{i_3 , j_3 }$
belongs to $\GG$.
Thanks to ($\sharp$) together with $g_1 \notin \GG$, 
we have $ (i_1, j_2) \in  S$.
Thanks to ($\sharp$) together with $g_2 \notin \GG$, we have $(i_1,j_1) ,(i_3,j_1)\in S$ and $(i_3,j_2) \notin S$.
Similarly, $g_3 \notin \GG$ implies 
$(i_2,j_1),(i_2,j_3) \in S$.
Thus,
$(i_1 ,j_2), (i_2 , j_3 )$, $(i_3 , j_1 )
\in S$
and $(i_3 , j_2) \notin S$.
This contradicts $f \in I_S$.

(2)
Suppose that none of
$g_1 = x_{i_1 , j_2 } x_{i_2 , j_1 } -  x_{i_1 , j_1 } x_{i_2 , j_2 }$,
$g_2 = x_{i_2 , j_2 } x_{i_3 , j_1 } -  x_{i_2 , j_1 } x_{i_3 , j_2 }$
and
$g_3 = x_{i_1 , j_3 } x_{i_2 , j_2 } -  x_{i_1 , j_2 } x_{i_2 , j_3 }$
belongs to $\GG$.
Thanks to ($\sharp$) together with $g_1 \notin \GG$, 
we have $ (i_2, j_2) \notin  S$.
Thanks to ($\sharp$) together with $g_2 \notin \GG$, we have $(i_2,j_1) \in S$ and $(i_3,j_1) \notin S$.
Similarly, $g_3 \notin \GG$ implies 
$(i_1,j_3) \notin S$.
Thus,
$(i_1 ,j_3), (i_2 , j_2 ),(i_3 , j_1 ) \notin S$
and $(i_2 , j_1) \in S$.
This contradicts $f \in I_S$.

(3)
Suppose that none of
$g_1 = x_{i_2 , j_3 } x_{i_3 , j_2 } -  x_{i_2 , j_2 } x_{i_3 , j_3 }$,
$g_2 = x_{i_1 , j_3 } x_{i_2 , j_2 } -  x_{i_1 , j_2 } x_{i_2 , j_3 }$
and
$g_3 = x_{i_2 , j_2 } x_{i_3 , j_1 } -  x_{i_2 , j_1 } x_{i_3 , j_2 }$
belongs to $\GG$.
Thanks to ($\sharp$) together with $g_1 \notin \GG$, 
we have $ (i_2, j_2) \in  S$.
Thanks to ($\sharp$) together with $g_2 \notin \GG$, we have $(i_1,j_3) \in S$ and $(i_2,j_3) \notin S$.
Similarly, $g_3 \notin \GG$ implies 
$(i_3,j_1) \in S$.
Thus,
$(i_1 ,j_3), (i_2 , j_2 ),(i_3 , j_1 ) \in S$
and $(i_2 , j_3 ) \notin S$.
This contradicts $f \in I_S$.

(4)
Suppose that none of
$g_1 = x_{i_1 , j_3 } x_{i_2 , j_1 } -  x_{i_1 , j_1 } x_{i_2 , j_3 }$
and
$g_2 = x_{i_1 , j_3 } x_{i_3 , j_2 } -  x_{i_1 , j_2 } x_{i_3 , j_3 }$
belongs to $\GG$.
Thanks to ($\sharp$) together with $g_1 \notin \GG$, 
we have $ (i_1, j_1) \in  S$.
Since $f \in I_S$, $\{(i_1,j_3),(i_2,j_1),(i_3,j_2)\} \cap S \neq \emptyset$.
Again, 
thanks to ($\sharp$) together with $g_1, g_2 \notin \GG$, we have $\{(i_1,j_3),(i_2,j_1),(i_3,j_2)\} \subset S$.
Since $f \in I_S$, we have $(i_3,j_3) \in S$.
This contradicts $g_2 \notin \GG$.

(5)
Suppose that none of
$g_1 = x_{i_1 , j_2 } x_{i_3 , j_1 } -  x_{i_1 , j_1 } x_{i_3 , j_2 }$
and
$g_2 = x_{i_2 , j_3 } x_{i_3 , j_1 } -  x_{i_2 , j_1 } x_{i_3 , j_3 }$
belongs to $\GG$.
Thanks to ($\sharp$) together with $g_1 \notin \GG$, 
we have $ (i_1, j_1) \in  S$.
Since $f \in I_S$, $\{(i_1,j_2),(i_2,j_3),(i_3,j_1)\} \cap S \neq \emptyset$.
Again, 
thanks to ($\sharp$) together with $g_1, g_2 \notin \GG$, we have $\{(i_1,j_2),(i_2,j_3),(i_3,j_1)\} \subset S$.
Since $f \in I_S$, we have $(i_3,j_3) \in S$.
This contradicts $g_2 \notin \GG$.

(6)
Suppose that none of
$g_1 = x_{i_1 , j_2 } x_{i_2 , j_1 } -  x_{i_1 , j_1 } x_{i_2 , j_2 }$
and
$g_2 = x_{i_2 , j_3 } x_{i_3 , j_2 } -  x_{i_2 , j_2 } x_{i_3 , j_3 }$
belongs to $\GG$.
Thanks to ($\sharp$) together with $g_2 \notin \GG$, 
we have $ (i_2, j_2) \in  S$.
This contradicts $g_1 \notin \GG$.

Thus, in all cases, $g_i \in \GG$ for some $i$.
Then $f$ is reduced to the cubic binomial $h$ by $g_i$
where monomials of $h$ have a common variable.
Hence $f$ is reduced to $0$ by $\GG$.

\medskip

On the other hand, suppose that $S$ is $2 \times 2$ block diagonal, that is, 
there exist integers $r$, $c$ such that 
$S=
\{ (i,j) \ | \ 1 \leq i \leq r, 1 \leq j \leq c \}
\cup 
\{ (i,j) \ | \ r < i \leq m, c <  j \leq n \}.
$
Let
$S'=
\{ (i,j) \ | \ 1 \leq i \leq r, 1 \leq j \leq c \}
$.
Then $S'$ is triangular.
Let
$
f=
x_{i,\ell}
x_{j,k}
-
x_{i,k}
x_{j,\ell}$
with
$ 
1\leq i < j \leq m,
1\leq k < \ell \leq n
$.
Then, $f \notin I_{S'}$
if and only if
$i \leq r < j$ and $k \leq c < \ell$
if and only if
$f \notin I_S$.
Hence we have $I_S=I_{S'}$.
\hfill{Q.E.D.}

\bigskip

\noindent
{\bf Acknowledgment.} \ 
The authors are grateful to
a referee for his useful suggestions
and
comments
for improving the expression of this paper.


\begin{thebibliography}{99}


\bibitem{CLO}
D. Cox, J. Little and D. O'Shea,
``Ideals, Varieties and Algorithms,'' Springer--Verlag,
Berlin, Heidelberg, New York, 1992.

\bibitem{HTY}
H. Hara, A. Takemura and R. Yoshida,
Markov bases for two-way subtable sum problems,
arXiv:math.CO/0708.2312v1, 2007.




\bibitem{quadratic}
H. Ohsugi and T. Hibi, 
Toric ideals generated by quadratic binomials,
{\em J. Algebra} {\bf 218} (1999), 509 -- 527. 




\bibitem{OHconti}
H. Ohsugi and T. Hibi,
Toric ideals arising from contingency tables,
{\it in} ``Commutative Algebra and Combinatorics,"
Ramanujan Mathematical Society Lecture Notes Series, Number 4,
Ramanujan Mathematical Society,
India, in press.



\bibitem{Stu}
B. Sturmfels, ``Gr\"obner bases and convex polytopes,"
Amer. Math. Soc., Providence, RI, 1995.


\end{thebibliography}
\end{document}